\documentclass[11pt]{article}

\usepackage[T1]{fontenc}
\usepackage{mathtools}
\usepackage{amssymb}
\usepackage{cite}
\usepackage{booktabs}
\usepackage{multirow}
\usepackage{caption, subcaption}
\usepackage{graphicx, float}
\usepackage{algorithm, algpseudocode}
\usepackage[margin=1in]{geometry}
\usepackage{titling}
\usepackage{fancyhdr}
\usepackage{makecell}

\usepackage{authblk}

\usepackage[colorlinks,urlcolor=blue,linkcolor=blue,citecolor=blue]{hyperref}

\title{\Large\bfseries Are Deep Learning Based Hybrid PDE Solvers Reliable?\\
Why Training Paradigms and Update Strategies Matter}

\author[1]{Yuhan Wu}
\author[2]{Jan Willem van Beek}
\author[2]{Victorita Dolean}
\author[1]{Alexander Heinlein}

\affil[1]{\small Delft Institute of Applied Mathematics, Delft University of Technology, Delft, The Netherlands}
\affil[2]{\small Department of Mathematics and Computer Science, Eindhoven University of Technology, The Netherlands}

\date{\vspace{-2em}}

\begin{document}

\maketitle

\thispagestyle{fancy}
\fancyhf{}

\fancyhead[C]{
  \parbox{0.9\textwidth}{\centering
  \fontsize{9}{4}\selectfont
  \copyright\ 2026 IEEE. Personal use of this material is permitted. Permission from IEEE must be obtained for all other uses, in any current or future media, including reprinting/republishing this material for advertising or promotional purposes, creating new collective works, for resale or redistribution to servers or lists, or reuse of any copyrighted component of this work in other works.
  }
}

\fancyfoot[C]{\thepage}

\renewcommand{\headrulewidth}{0pt}

\begin{abstract}
    Deep learning-based hybrid iterative methods (DL-HIMs) integrate stationary iterative solvers with neural operators, utilizing their complementary spectral biases to accelerate convergence. However, many DL-HIMs stagnate at false fixed points where updates vanish while significant physical residuals persist. We demonstrate that reliability depends on training paradigms and update strategies rather than neural architectures alone. Through a detailed study of a nonlinear DeepONet-based and a linear FFT-based DL-HIM, we show that convergence degrades when training objectives are misaligned with solver dynamics and problem physics. We examine Anderson acceleration and find its classical form ill-suited for nonlinear neural operators. We introduce physics-aware Anderson acceleration (PA-AA), which minimizes the physical residual rather than the fixed-point update. Numerical experiments confirm that PA-AA restores reliable convergence in substantially fewer iterations. These findings suggest that addressing controversies surrounding AI-based PDE solvers requires attention not only to architectures but to physically informed training and iteration design.
\end{abstract}

\section{Introduction}
Solving large-scale linear systems resulting from the discretization of partial differential equations (PDEs) robustly and efficiently is essential for scientific and engineering applications.

Stationary iterative methods for solving the linear system $Au=f$ takes the form:
$$
    u^{(k+1)} = u^{(k)} + \omega B(f - Au^{(k)})
$$
where $\omega$ denotes the relaxation parameter and $B$ is a cheap-to-apply preconditioner; cf.~\cite{van2003iterative}. The convergence is determined by the spectral properties of the error propagation operator $I - \omega BA$.

Preconditioners aim to approximate the inverse of the system matrix, $B\approx A^{-1}$. Classical matrix splitting preconditioners, such as Jacobi and Gauss-Seidel~\cite{van2003iterative}, are inexpensive but converge slowly for low-frequency error components and may even diverge on nonsymmetric or indefinite systems~\cite{van2003iterative}. More sophisticated preconditioners include incomplete factorizations~\cite{benzi2002preconditioning}, multigrid~\cite{trottenberg_multigrid_2000}, and domain decomposition~\cite{dolean2015introduction}.

Beyond preconditioning techniques, the relaxation parameter can be optimized to accelerate convergence; see, e.g., the Chebyshev semi-iterative approach~\cite{golub1961chebyshev}. Alternatively, Krylov subspace methods~\cite{van2003iterative} can significantly improve convergence. Similarly, multisecant information from recent iterates for general (potentially nonlinear) fixed-point iteration schemes can be leveraged; cf., e.g., Anderson acceleration (AA)~\cite{walker2011anderson}. For complex problems, convergence remains sensitive to the spectral properties of the (preconditioned) system matrix, so effective preconditioners remain critical.

Additionally, deep learning has gained considerable attention for solving differential equations, including physics-informed neural networks (PINNs)~\cite{raissi2019physics} and neural operators such as Fourier neural operators (FNOs)~\cite{li2020fourier} and deep operator networks (DeepONets)~\cite{lu2021deeponet}. However, deep neural networks typically exhibit a spectral bias~\cite{rahaman_spectral_2019}: low-frequency modes converge faster than high-frequency modes during training. Since classical stationary smoothers (e.g., damped Jacobi and Gauss-Seidel) primarily reduce high-frequency errors, this complementarity motivates deep learning-based hybrid iterative methods (DL-HIMs), which alternate a classical stationary smoother with a neural operator.

This idea was first formalized as the hybrid iterative numerical transferable solver (HINTS) in Zhang et al.~\cite{zhang2024blending}, which alternates $n$ steps of a stationary smoother with a DeepONet model. Hu and Jin~\cite{hu2025hybrid} later replaced DeepONet with a multiple-input operator network (MIONet)~\cite{jin2022mionet}, as MIONet preserves linearity with respect to the right-hand side $f$. Similarly, Sun et al.~\cite{sun2025learning} leverage learned singularity-encoded Green's functions as the neural component in the DL-HIM. To accurately target the error components unresolved by the smoother, Cui et al.~\cite{cui2025hybrid} introduced the Fourier neural solver (FNS), which utilizes a specialized fast Fourier transform (FFT)-based architecture and end-to-end training, albeit currently limited to structured grids.

Despite these advances, a fundamental challenge persists: the spectral gap between the classical smoother and the neural operator often leads to convergence stagnation \cite{kopanivcakova2025deeponet}. Considering the DL-HIM as a fixed point iteration $u^{(k+1)} = G(u^{(k)})$, this gap may manifest as 'false fixed points' where updates from both the neural networks and the classical smoother vanish while the actual residual remains large.

Existing research emphasizes neural architectures, whereas training paradigms and update strategies remain under-explored. In DeepONet/MIONet-based approaches, the neural operator is typically trained offline on static datasets, which often require pre-solved PDE solutions and introduce a distribution mismatch between training and inference. By contrast, FNS uses an end-to-end dynamic training framework and residual-based objectives that avoid the pre-solving cost and train-inference mismatch, but typically increase computational and memory overhead in training.

Moreover, a full-cycle DL-HIM update defines a fixed-point iteration in which each cycle applies the same composition of smoother and neural correction, $u^{(k+1)}=G(u^{(k)})$. While the underlying PDE may be linear, the neural operator may still exhibit nonlinearity with respect to the residual. This motivates the use of acceleration strategies such as AA that apply to general nonlinear fixed-point iterations.

This work provides a systematic study of training paradigms and update strategies for two representative DL-HIM solvers, HINTS and FNS. The main contributions are as follows:
\begin{itemize}
    \item We characterize the DL-HIM iteration as a fixed point iteration $u^{(k+1)}=G_\theta(u^{(k)})$, which reveals that the inherent spectral gap leads to a misalignment between the mathematical fixed point of the iterator and the true solution of the discretized PDE. This provides a mechanistic understanding of convergence stagnation, identifying the false fixed point as one of the primary barriers in DL-HIMs.

    \item While recent research has primarily focused on designing neural architectures, the impact of training objectives remains under-explored. We evaluated different training paradigms, revealing that the effectiveness of a training objective is highly context-dependent. Our experiments demonstrate that DeepONet-based DL-HIM solvers tend to favor residual-based losses, whereas FNS may benefit from error-based objectives.

    \item We analyze the impact of the update strategies on DL-HIMs, including exact line search and AA. We highlight that the nonlinearity of neural operators hinders standard AA from effectively escaping false fixed points. To address this, we propose physics-aware AA (PA-AA), a strategy that explicitly minimizes the physical residual norm $\|f-Au\|$ rather than the fixed-point update. In our numerical experiments, PA-AA breaks the stagnation barrier for nonlinear DL-HIM solvers, enabling them to converge stably to a much smaller physical residual with fewer iterations.

\end{itemize}

The code for this work can be found in the GitHub repository \url{https://github.com/HNU-WYH/DL-HIM}.

\section{FRAMEWORK OF DL-HIM}\label{sec:framework}
\subsection{Problem Setup}
Consider the following stationary linear parametric PDE:
\begin{equation}
    \mathcal{L}\left[k(x), u(x)\right] = f(x)
    \label{eq:PDE}
\end{equation}
where $\mathcal L$ denotes the differential operator, $u(x)$ is the solution, $k(x)$ is the parameter function, and $f(x)$ is the source term. We assume that $k$ and $f$ are such that the problem is well posed and admits a unique solution $u$.

After discretization, Eq.~\eqref{eq:PDE} yields a system of linear equations:
\begin{equation}
    Au = f
    \label{eq:LA}
\end{equation}
where $A\in\mathbb R^{n\times n}$ is the discretized linear operator depending on the parameter function $k(x)$. While Eq.~\eqref{eq:LA} is linear, the hybrid solver itself may incorporate nonlinear components, leading to a complex interplay between the linear system and the nonlinear hybrid solver.

\subsection{Fixed-Point Iteration in DL-HIM}
A standard full DL-HIM cycle alternates $n$ steps of a classical smoother with one step of a neural operator correction. For a given iterate $u^{(i)}$, the smoother and neural updates are:
\begin{equation}\label{eq:dl-him}
    \begin{cases}
        \tilde u = u^{(i)} + \omega \mathcal{S}(f - A u^{(i)}) \      & {\text{(repeat $n$ times)}} \\
        u^{(i+1)} = \tilde u + \mathcal{N_\theta} (f - A \tilde u) \  & {\text{(neural operator)}}
    \end{cases}
\end{equation}
where $\mathcal S$ is the preconditioner for the stationary smoother (e.g., Jacobi or Gauss-Seidel), $\mathcal{N}_\theta$ is the neural operator parameterized by $\theta$, and $\omega$ is the relaxation parameter for the smoother. The neural operator does not require an explicit relaxation parameter, as its output scale is implicitly learned during training and can be interpreted as an adaptive correction.

The solving process in Eq.~\eqref{eq:dl-him} can be defined as a fixed point mapping $G_\theta: \mathbb R^n\to \mathbb R^n$ such that:
\begin{equation}
    u^{(k+1)} = G_\theta(u^{(k)}) = \mathcal{M}_{\mathcal{N}_\theta} \circ \mathcal{M}_{\mathcal{S}}^n (u^{(k)})
    \label{eq:fixed_point_map}
\end{equation}
where $\mathcal{M}_{\mathcal{S}}$ and $\mathcal{M}_{\mathcal{N}_\theta}$ represent the mapping induced by one step of a stationary smoother $\mathcal S$ and a neural operator $\mathcal{N}_\theta$, separately.

Let $\delta_k = G_\theta(u^{(k)}) - u^{(k)}$ represent the update of a full DL-HIM cycle. The DL-HIM iteration terminates when the update $\delta^{(k)}$ vanishes, signaling that a mathematical fixed point has been reached. Ideally, the mathematical fixed point of $G_\theta$ should coincide with the PDE solution such that
$$
    \delta^{(k)} = 0 \iff r^{(k)} = 0.
$$

However, this equivalence may break down in DL-HIMs due to the inherent spectral gap. On fine meshes, stationary smoothers produce negligible updates for low-frequency errors. Simultaneously, the neural operator may also fail to respond to certain error modes, yielding near-zero corrections or unreliable updates, which are subsequently eliminated by the stationary smoother. Ultimately, this results in a negligible full-cycle update $\delta^{(k)}$ in DL-HIMs.

As a result, it is possible that the iteration traps at a false fixed point where $\delta^{(k)} \approx 0$ despite a significant physical residual $r^{(k)}$, meaning that the fixed point iteration $G_\theta$ fails to satisfy the underlying PDE. We show an example of this phenomenon occurring in Section~ \nameref{sec:numerical_experiments}. Moreover, we develop a variant of AA in Section~ \nameref{sec:update_strategy} to address this problem.

\subsection{Error and Residual Propagation Operators}
To analyze its convergence, we can approach this iteration from the perspectives of error propagation and residual propagation.

\subsubsection{Error Propagation}
Let $u^\star$ denote the exact solution of Eq.~\eqref{eq:LA}, and $e^{(k)} = u^\star - u^{(k)}$ be the error at the $k$-th cycle. The evolution of the error through a full DL-HIM cycle is governed by the following components:

\begin{itemize}
    \item \textbf{Stationary smoother}: For $n$ steps of the smoother $\mathcal{S}$, the error propagation operator $E_{\mathcal{S}}^n$ (mapping the current error to the next-iteration error induced by the given update rule) is given by
          \begin{equation}
              E_{\mathcal{S}}^n \coloneqq (I - \omega \mathcal{S}A)^n
          \end{equation}

    \item \textbf{Neural operator}: the error propagation through the neural operator step defines a mapping $E_{\mathcal{N}_\theta}$ that acts on the error $e$ as
          \begin{equation}
              E_{\mathcal{N}_\theta}(e) \coloneqq e - \mathcal{N}_\theta(A e);
          \end{equation}
\end{itemize}

The full-cycle error propagation operator $E$ for one full DL-HIM cycle is defined by the composition:
\begin{equation}
    e^{(k+1)} = E(e^{(k)}) \coloneqq (E_{\mathcal{N}_\theta} \circ E_{\mathcal{S}}^n)(e^{(k)})
\end{equation}

If the architecture of the neural operator is nonlinear (e.g., DeepONet), the resulting full-cycle error propagation operator $E$ is also nonlinear. In that case, the convergence rate may be state-dependent, rather than being determined by a global spectral radius.

\subsubsection{Residual Propagation}
Let $r^{(k)} = f - A u^{(k)}$ be the residual during the $k$-th cycle. Given the relationship $r^{(k)} = Ae^{(k)}$, the residual propagation follows a similar structure to that of the error:
\begin{itemize}
    \item \textbf{Stationary smoother}:
          \begin{equation}
              R_S^n \coloneqq (I - \omega A \mathcal S)^n
          \end{equation}

    \item \textbf{Neural operator}:
          \begin{equation}
              R_{\mathcal{N}_\theta}(r) \coloneqq r - A \mathcal{N}_\theta(r)
          \end{equation}
\end{itemize}

The full-cycle residual propagation operator $R$ is thus defined as the composition:
\begin{equation}
    r^{(k+1)} = R(r^{(k)}) \coloneqq (R_{\mathcal{N}_\theta} \circ R_S^n)(r^{(k)})
\end{equation}

Specifically, the relationship between the residual propagation operator $R$ and the error propagation operator $E$ is given by:
$$
    R = A\circ E\circ A^{-1}.
$$
While $E$ and $R$ share the same spectral properties in the linear case, for nonlinear DL-HIMs, the simple spectral equivalence no longer holds.

\section{TRAINING PARADIGMS OF DL-HIM}

Ideally, the neural operator parameters $\theta$ should be chosen such that each DL-HIM cycle consistently reduces the error or residual. The objective can be formalized as minimizing the contraction factors $q_r$ or $q_e$:
\begin{equation}
    \begin{aligned}
        \|Rr^{(k)}\|_Y & \le\,q_r\|r^{(k)}\|_Y \\
        \|Ee^{(k)}\|_X & \le\,q_e\|e^{(k)}\|_X
    \end{aligned}
    \label{eq:contraction-residual}
    \quad \text{where } 0 \le q_r, q_e < 1
\end{equation}
where $\|\cdot\|_X$ and $\|\cdot\|_Y$ denote the chosen norms in the solution and residual spaces, respectively.

For linear neural operators (e.g., FNS and MIONet), this is equivalent to minimizing the spectral radius $\rho(E) < 1$ or $\rho(R) < 1$. However, directly optimizing the spectral radius is computationally intractable, as it depends implicitly on the full operator composition and requires repeated eigenvalue estimations of large-scale operators. Therefore, tractable loss functions are required as proxies for these contraction factors.

\subsection{Static Training Framework}\label{sec:static_training}
In the static framework, $\mathcal N_\theta$ is trained offline as a single-step solver on a precomputed dataset $\{(A_j, f_j, u_j^\star)\}_{j=1}^N$ with $u_j^\star = A_j^{-1}f_j$. To minimize the cost of generating reference solutions, the training is typically performed on coarse grids. The resolution-invariant nature of neural operators enables architectures like DeepONet and FNS to be subsequently deployed on finer grids during inference.

A key limitation of static training is the distribution mismatch: the model is only trained on initial residuals $r^{(0)}=f$ but acts on evolved residuals $r^{(k)}$ produced by the hybrid iteration during inference. This mismatch arises during the iterative process: the spectral gap between the smoother and the neural operator dictates which error components stagnate. The persistence of these unresolved modes drives the continuous evolution of the residual, steadily shifting its distribution away from the initial training data.

Since the static framework cannot expose the model to the true inference-stage distribution, the choice of the training objective becomes critical. We consider two primary families of loss functions in our assessment:

\paragraph{Error-Based Loss} This loss penalizes the distance to the reference solution. Given a dataset $\{(A_j,f_j,u_j^\star)\}_{j=1}^N$, a standard choice (e.g., HINTS) is the relative solution loss:
\begin{equation}
    \begin{aligned}
        \mathcal L_{\mathrm{err}}^\mathrm{stat}(\theta) & = \frac1N \sum_{j=1}^N \frac{\bigl\|u_j^\star - \mathcal N_\theta(f_j)\bigr\|_X^2} {\|u_j^\star\|_X^2}
        \\& \approx \mathbb E \left[ \frac{\bigl\|E_{\mathcal N_\theta} \left(e_j^{(0)}\right)\|_X^2} {\|e_j^{(0)}\|_X^2}\right]
    \end{aligned}
\end{equation}
where $\|\cdot\|_X$ denotes a chosen norm, and $E_{\mathcal N_\theta}(e)=e-\mathcal N_\theta(Ae)$ is the error propagation operator. The equality follows from $f_j= Au_j^\star=Ae^{(0)}$ and the definition of $E_{\mathcal N_\theta}$. Minimizing this loss optimizes $E_{\mathcal N_\theta}$ by training $\mathcal N_\theta$ to reduce error modes in the solution subspace spanned by the training data.

\paragraph{Residual-Based Loss} This loss penalizes the violation of the PDE directly, eliminating the need for reference solutions. Given training samples $\{(A_j,f_j)\}_{j=1}^N$, the relative residual loss takes the form:
\begin{equation}
    \begin{aligned}
        \mathcal L_{\mathrm{res}}^\mathrm{stat}(\theta) & = \frac1N \sum_{j=1}^N \frac{\bigl\|f_j - A_j \mathcal N_\theta(f_j)\bigr\|_Y^2}{\|f_j\|_Y^2}
        \\& \approx \mathbb{E} \left[\frac{\bigl\|R_{\mathcal N_\theta}\left(r_j^{(0)}\right)\bigr\|_Y^2}{\|r_j^{(0)}\bigr\|_Y^2}\right]
    \end{aligned}
\end{equation}
where $f_j=r^{(0)}$ is the initial residual with a zero guess, $R_{\mathcal N_\theta}(r) = r - A\mathcal N_\theta(r)$ is the residual propagation operator for $\mathcal N_\theta$, and $\|\cdot\|_Y$ is a norm on the space of residuals. This loss serves as an empirical proxy for the contraction factor of the residual propagation.

The choice of norms in the loss further influences which spectral modes are emphasized during training. In our experiments, we consider the following choices:
\begin{itemize}
    \item $\ell^2$ norms on solutions or residuals;
    \item $\ell^1$ norms on solutions or residuals;
    \item $H^1$-type norms of the form
          $$
              \|x\|_{H^1}^2 = \|x\|_2^2 + \lambda \|\nabla_h x\|_2^2
          $$
          where $\nabla_h$ denotes the discrete gradient operator (implemented via central differences) and $\lambda$ is a hyperparameter.
\end{itemize}

In the Fourier domain, an $H^1$-type norm weights each mode $\xi$ by $(1+\lambda |\xi|^2)$, thereby increasingly penalizing highly oscillatory components. Similarly, a residual-based $\ell^2$ loss $\|r\|^2_2$ is equivalent to an $A^\top A$-induced energy norm on the error, $\|e\|^2_{A^\top A} = \|Ae\|^2_2$. For elliptic PDEs, such as the Poisson equation, it also amplifies high-frequency components more strongly than low-frequency ones, thus biasing training toward reducing highly oscillatory errors.

In summary, error-based and residual-based losses provide two complementary ways of shaping the single-step error and residual propagation operators of $\mathcal N_\theta$ on the subspaces spanned by the training data, while different norm choices act as frequency weighting mechanisms within these subspaces. We empirically compare these design choices in Section~ \nameref{sec:experiment_loss_functions}.

\subsection{Dynamic Training Framework}\label{sec:dynamic_training}
The dynamic end-to-end framework explicitly unrolls $K$ DL-HIM cycles within a differentiable computational graph. This shifts the training objective from single-step accuracy of a neural operator to optimizing K-step convergence for the full-cycle propagation operator $E$ or $R$.

By simulating the actual inference workflow, the dynamic framework naturally exposes $\mathcal N_\theta$ to the evolving residuals encountered over multiple steps, thereby mitigating differences in the distribution of inputs between training and inference. This also ensures that the neural operator targets error components unresolved by the classical smoother, thereby tailoring $\mathcal N_\theta$ to correct those difficult spectral modes.

In each training episode, we sample a PDE instance and a right-hand side (RHS) $f$, start from $u^{(0)}=0$, and unroll the solver for $K$ cycles so that $\mathcal N_\theta$ receives residuals $\{r^{(k)}\}_{k=0}^{K-1}$ generated on-the-fly. Accordingly, we typically define cumulative losses over the $K$ DL-HIM cycles, aggregating per-step losses along the unrolled trajectory. These cumulative losses serve as a proxy for the multi-step contraction factor.

However, unrolling $K$ iterations incurs significantly higher memory consumption and computational overhead compared to static training. The elongated computational graph may also exacerbate training instability, such as gradient explosion or vanishing. Moreover, the generalization challenge may still exist when the finite horizon of $K$ steps is not sufficiently large.

\paragraph{Error-Based Loss}
When reference solutions $u_i^\star$ are available, one can directly minimize the accumulated normalized solution error over $K$ iterations
\begin{equation}
    \label{eq:dyn_loss_sol}
    \begin{aligned}
        \mathcal{L}_{\mathrm{sol}}^\mathrm{dyn}(\theta) & =\frac{1}{NK}\sum_{i=1}^{N}\sum_{k=1}^{K}\frac{\|u_{i}^{\star}-u_{i}^{(k)}\|_{X}^{2}}{\|u_{i}^{\star}\|_{X}^{2}}
        \\&=\frac{1}{NK}\sum_{i=1}^{N}\sum_{k=1}^{K}\frac{\|E^{k}e_i^{(0)}\|_{X}^{2}}{\|e_{i}^{(0)}\|_{X}^{2}},
    \end{aligned}
\end{equation}
where $u^{(k)}_i$ is the intermediate solution of the $i$-th training problem after $k$ DL-HIM cycles, $e_i^{(k)} = u_i^\star - u_i^{(k)}$ is the corresponding error, and $E$ denotes the full-cycle error propagation operator. While this loss directly measures the average error contraction across multiple iterations, it requires pre-solved reference solutions $u_i^\star$.

\paragraph{Residual-Based Loss}
To avoid the cost of pre-solving PDEs, the residual-based dynamic loss optimizes the full-cycle residual propagation operator directly and takes the form:
\begin{equation}
    \begin{aligned}
        \mathcal{L}_{\mathrm{res}}^\mathrm{dyn}(\theta) & =\frac{1}{NK}\sum_{i=1}^{N}\sum_{k=1}^{K}\frac{\|r_{i}^{(k)}\|_{Y}^{2}}{\|f_{i}\|_{Y}^{2}}
        \\&=\frac{1}{NK}\sum_{i=1}^{N}\sum_{k=1}^{K}\frac{\|R^{k}r_i^{(0)}\|_{Y}^{2}}{\|r_{i}^{(0)}\|_{Y}^{2}},
    \end{aligned}
\end{equation}
This loss depends only on $(A_i,f_i)$ without $u_i^\star$. Since $r_i^{(k)} = A_i e_i^{(k)}$, $\mathcal{L}_{\mathrm{res}}^\mathrm{dyn}$ is equivalent to an $A^\top A$-weighted error loss, which implicitly introduces a bias towards reducing oscillatory error components.

\section{UPDATE STRATEGIES OF DL-HIM}\label{sec:update_strategy}
In addition to the training paradigm and architecture of neural operators, the actual performance of hybrid solvers also critically depends on the update strategy for moving the current iterate $u^{(k)}$ to $u^{(k+1)}$ during inference. This has remained relatively under-explored in most existing studies. In this section, we investigate the impact of update strategies on the convergence of DL-HIM, focusing on step size selection and multi-secant acceleration techniques to mitigate convergence stagnation.

\subsection{Step Size Selection}
The neural operator $\mathcal{N}_\theta$ is explicitly trained to approximate the error corresponding to the given current residual. Consequently, it is natural to apply this correction directly without scaling, and a fixed step size $\alpha_k=1$ is widely adopted in DL-HIMs such as HINTS\cite{zhang2024blending}.

In contrast, FNS\cite{cui2025hybrid} introduces a scalar step size $\alpha_k$ to modulate the update magnitude along the direction provided by the neural operator $p^{(k)} = \mathcal{N}_\theta(r^{(k)})$.

The selection of $\alpha_k$ depends on the properties of the system matrix $A$. For symmetric positive definite (SPD) matrices, FNS adopts an exact line search strategy derived from the conjugate gradient (CG) method, which determines the optimal $\alpha_k$ by minimizing the $A$-norm of the error $e_k$ along the direction $p_k$:
\begin{equation}\label{eq:adaptive_step}
    \alpha_k = \arg \min_{\alpha} \|e_{k} - \alpha p_k\|^2_A = \frac{p_k^\top r_k}{p_k^\top A p_k}
\end{equation}

This adaptive step size exploits the update direction provided by $\mathcal{N}_\theta$ more effectively. However, for non-symmetric or indefinite problems (e.g., the Helmholtz equation), the energy norm $\|\cdot\|_A$ is ill-defined, and the denominator in Eq.~\eqref{eq:adaptive_step} may vanish or become negative. We therefore treat its use in non-SPD settings as heuristic and evaluate its behavior empirically in Section~\nameref{sec:experiment_update}.

\subsection{Anderson Acceleration (AA)}
While adaptive step sizes can optimize the update magnitude, their reliance on the $A$-norm limits their applicability to SPD problems. Moreover, scalar tuning cannot correct the update direction if the neural operator itself provides an incorrect prediction. This motivates the use of multisecant acceleration methods, which leverage historical information to construct better update directions without assuming specific properties of the operator $A$.

Anderson acceleration (AA)~\cite{walker2011anderson} is a mature numerical technique designed to accelerate the convergence of fixed-point iterations. AA leverages multisecant information by storing the history of the residuals corresponding to the $m$ most recent iterates. An accelerated iterate is then given by an optimal linear combination of this history.

A standard implementation of AA applied to the DL-HIM solver $G_\theta$ is detailed in Algorithm~\ref{alg:aa}.

\begin{algorithm}[h!]
    \caption{Standard AA for DL-HIM}
    \label{alg:aa}
    \begin{algorithmic}[1]
        \Require Initial guess $u_0$, memory size $m \ge 1$
        \State $u_1 \gets G_\theta(u_0)$
        \For{$k = 1, 2, \dots$ \textbf{until convergence}}
        \State $g_k \gets G_\theta(u_k)$
        \State $r_k \gets g_k - u_k$    \Comment{Fixed-Point Residual}
        \State $m_k \gets \min(m, k)$

        \State $\mathcal{R}_k \gets [r_k, r_{k-1}, \dots, r_{k-m_k}]$
        \State $\mathcal U_k \gets [u_k - u_{k-1}, \dots, u_k - u_{k-m_k}]$

        \State find $\alpha = (\alpha_0, \dots, \alpha_{m_k})^T$ that minimizes:
        $$
            \quad \left\| \sum_{j=0}^{m_k} \alpha_j r_{k-j} \right\|_2 \
            \text{subject to } \sum_{j=0}^{m_k} \alpha_j = 1
        $$

        \State $u_{k+1} \gets \sum_{j=0}^{m_k} \alpha_j g_{k-j}$
        \EndFor
    \end{algorithmic}
\end{algorithm}

\subsection{Physics-Aware Anderson Acceleration (PA-AA)}
Although standard AA can be applied to general nonlinear fixed point mappings, it is designed to minimize the norm of $\delta^{(k)}$ and may inadvertently accelerate the nonlinear DL-HIM solver toward a false fixed point such that the neural operator's corrections have ceased, even though the underlying PDE is not satisfied.

To solve this problem, we proposed physics-aware AA (PA-AA). The core idea is to directly minimize the residual of the linear system rather than the update. The implementation of PA-AA is summarized in Algorithm~\ref{alg:pa_aa}.
\begin{algorithm}[h!]
    \caption{PA-AA for DL-HIM}\label{alg:pa_aa}
    \begin{algorithmic}[1]
        \Require Initial guess $u_0$, memory size $m \ge 1$, operator $A$, RHS $f$
        \State $u_1 \gets G_\theta(u_0)$
        \For{$k = 1, 2, \dots$}
        \State $g_k \gets G_\theta(u_k)$
        \State $r_k \gets f - A g_k$    \Comment{Physical Residual}
        \State $m_k \gets \min(m, k)$
        \State $\mathcal{R}_k \gets [r_k, r_{k-1},\dots, r_{k-m_k}]$
        \State $\mathcal{G}_k \gets [g_k, g_{k-1}, \dots, g_{k-m_k}]$

        \State find $\alpha = (\alpha_0, \dots, \alpha_{m_k})^T$ that minimizes:
        $$
            \quad \left\| \sum_{j=0}^{m_k} \alpha_j r_{k-j} \right\|_2 \
            \text{subject to } \sum_{j=0}^{m_k} \alpha_j = 1
        $$
        \State $u_{k+1} \gets \sum_{j=0}^{m_k} \alpha_j g_{k-j}$
        \EndFor
    \end{algorithmic}
\end{algorithm}

\section{NUMERICAL EXPERIMENTS} \label{sec:numerical_experiments}
In this section, we provide a critical assessment of the proposed frameworks. We evaluate performance on two representative parametric PDEs with homogeneous Dirichlet boundary conditions on $\Omega=[0,1]$:
\begin{enumerate}
    \item Stochastic diffusion equation:
          \begin{equation}
              -\nabla\cdot\left(k(x)\nabla u(x)\right)  =f(x)
          \end{equation}
          where $k(x)>0$ is a random coefficient field that yields a symmetric positive definite (SPD) system.

    \item Indefinite Helmholtz equation:
          \begin{equation}
              -\Delta u(x) - k^2(x)u(x) = f(x)
          \end{equation}
          where $k(x)$ represents the wavenumber function. We select sufficiently large wavenumbers to ensure that the discretized system is challenging for standard solvers.
\end{enumerate}

\subsection{Data Generation and Setup}
To construct the training and testing datasets, we sample the parameter function $k(x)$ and the source term $f(x)$ from Gaussian random fields (GRFs) with a radial basis function (RBF) covariance
\begin{equation}
    \text{Cov}({x}_1, {x}_2) = \sigma^2 \exp\left( - \frac{\|x_1 - x_2\|^2}{2l^2} \right),
\end{equation}
where $l$ denotes the correlation length and $\sigma$ denotes the standard deviation. To ensure the physical validity (e.g., ellipticity for the diffusion equation), we shift and clip $k(x)$ as
\begin{equation}
    k(x) = \max({\text{GRF}}(x) + \mu_k, k_{\min}),
\end{equation}
while $f(x)$ is sampled with zero mean and no truncation. We use $(\mu_k,k_{\min},\sigma_k,l_k)=(1.0,0.3,0.3,0.1)$ for diffusion and $(8.0,3.0,2.0,0.2)$ for Helmholtz; for the source term $f(x)$, we set $(\sigma_f,l_f)=(1.0,0.1)$ for both problems.

Neural operators are trained on a coarse grid ($N_{train}=31$) to minimize offline costs and are evaluated on fine grids of different sizes in numerical experiments. We employ a zero initial guess and empirically alternate $n=19$ Jacobi or Gauss-Seidel steps with one neural correction per DL-HIM cycle. The convergence of DL-HIM solver is measured by the norms of residual $\|r^{(k)}\|_2$  and error $\|e^{(k)}\|_2$.

\subsection{The False Fixed Point}
A primary insight of this work is the identification of false fixed points in the training of DL-HIMs. To illustrate this, we analyze the convergence performance of standard HINTS on the 1D Diffusion equation with a high-resolution grid of $N = 801$.
\begin{figure}[H]
    \centering
    \includegraphics[width=0.6\linewidth]{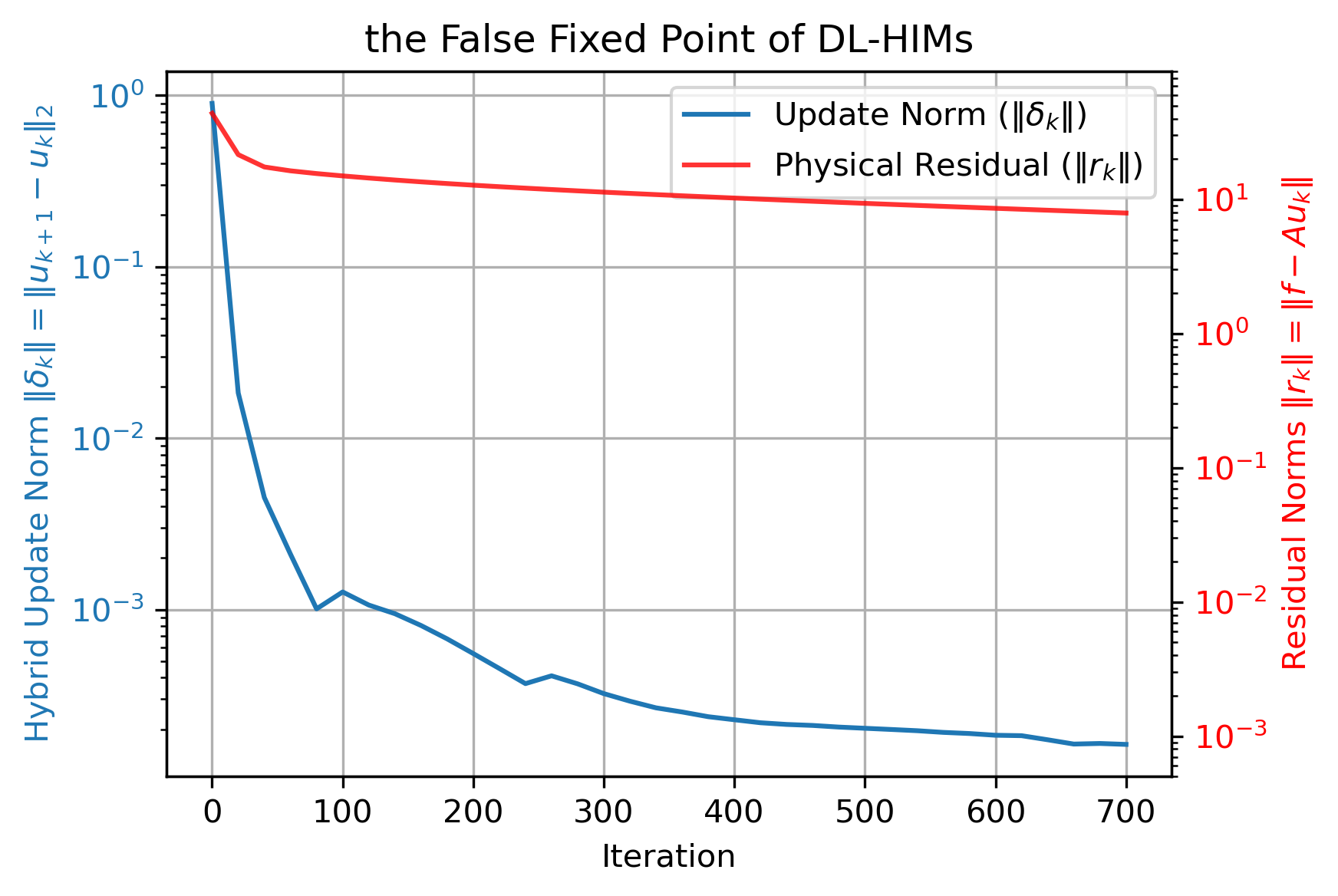}
    \caption{Update and Residual Norms of HINTS on the 1D Diffusion Equation with $N=801$}
    \label{fig:false_fixed_point}
\end{figure}

Figure~\ref{fig:false_fixed_point} depicts the evolution of the Update Norm $\|\delta_k\| = \|u_{k+1} - u_k\|_2$ and the Physical Residual $\|r_k\| = \|f - Au_k\|_2$ during the DL-HIM iterations. The update norm drops sharply and plateaus at approximately $10^{-4}$, while the physical residual stagnates at a high magnitude of approximately $10^{1}$ after an initial decrease.

The stagnation of both residual and update illustrates the existence of a false fixed point in the training of DL-HIMs. Especially on high-resolution grids ($N=801$), the spectral gap between the classical smoother and the neural operator becomes pronounced, with the classical smoother exhibiting slow convergence and the neural operator failing to handle some error modes unresolved by the smoother, leading to near-zero corrections even when the physical error is still significant.

\subsection{Effect of Loss Functions}\label{sec:experiment_loss_functions}
This section investigates how the training objective shapes the performance of the neural operator within the DL-HIM framework. To isolate the effect of the objective function from the training framework, we restrict the experiment to a static training framework. Our experiments reveal that the optimal loss function is highly architecture-dependent.

\begin{figure}[htbp]
    \label{fig:loss_hints_comparison}
    \centering
    \begin{subfigure}{0.48\linewidth}
        \centering
        \includegraphics[width=\textwidth]{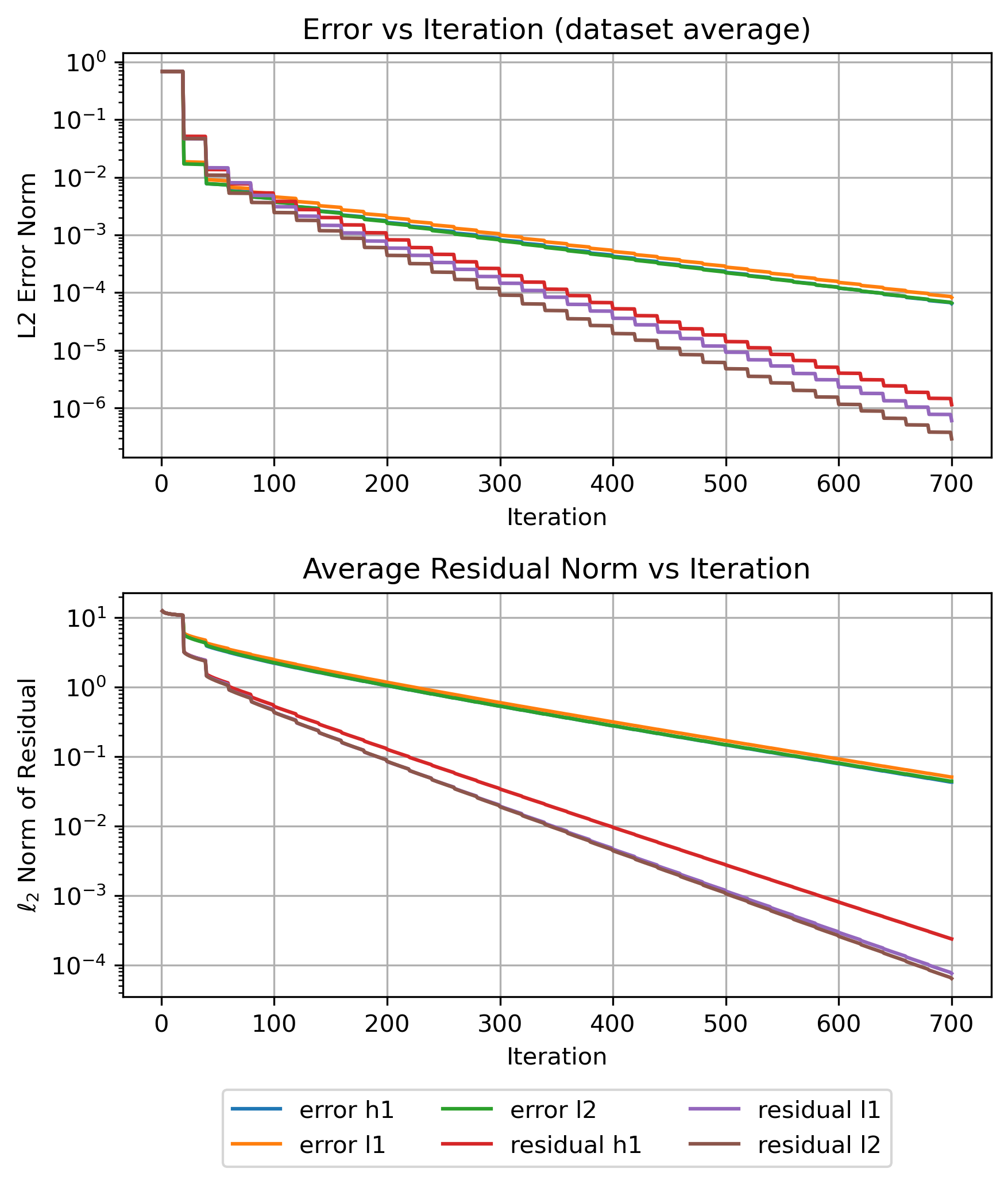}
        \caption{Convergence Behavior of HINTS on Diffusion Equations ($N=201$) with DeepONet Models Statically Trained by Different Loss Functions.}
        \label{fig:loss_hints_diffusion}
    \end{subfigure}
    \hfill
    \begin{subfigure}{0.48\linewidth}
        \centering
        \includegraphics[width=\textwidth]{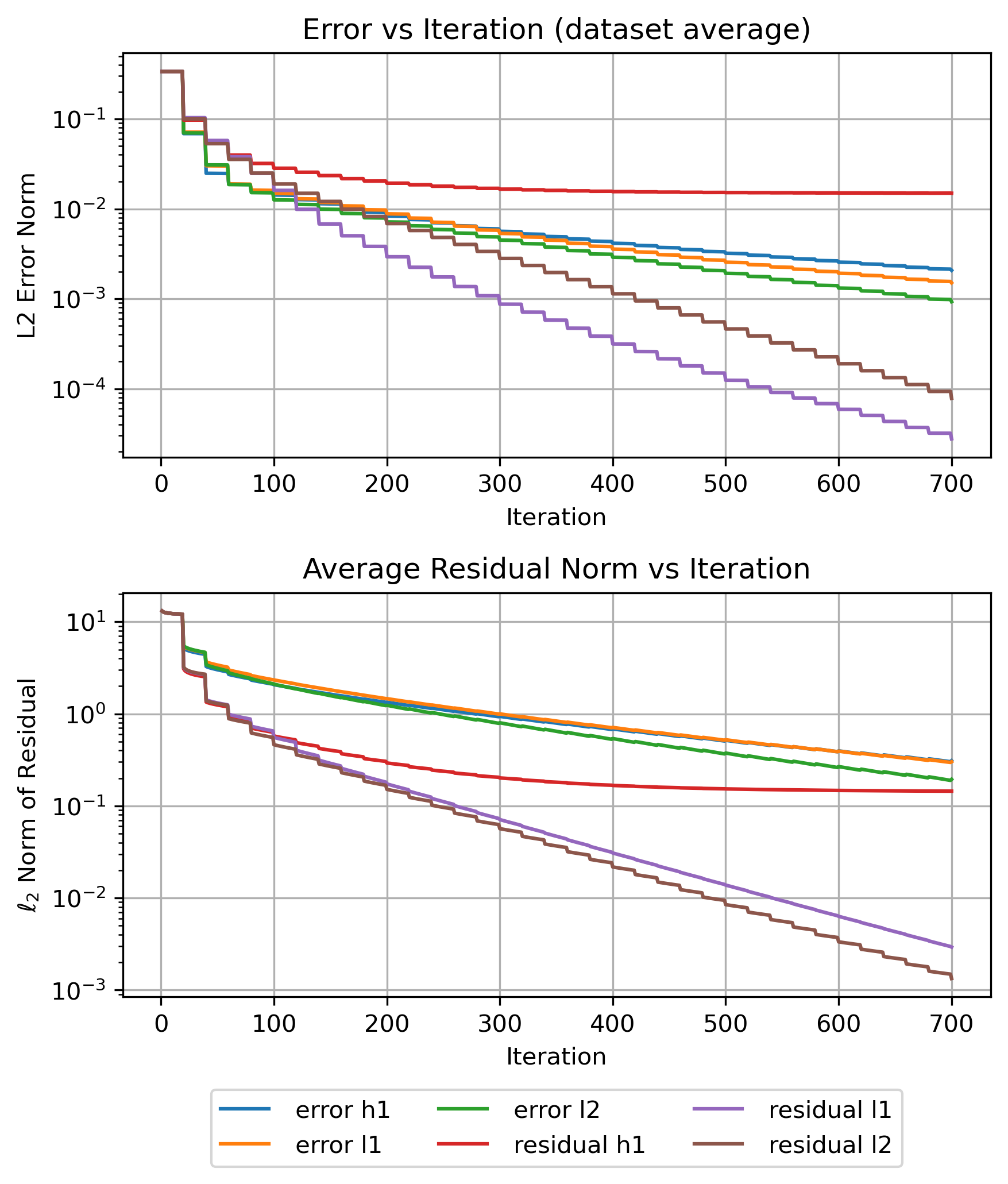}
        \caption{Convergence Behavior of HINTS on Helmholtz Equations ($N=201$) with DeepONet Models Statically Trained by Different Loss Functions.}
        \label{fig:loss_hints_helmholtz}
    \end{subfigure}
    \caption{Convergence Behavior of HINTS in the Static Training Framework}
\end{figure}

Figures~\ref{fig:loss_hints_diffusion} and \ref{fig:loss_hints_helmholtz} compare the convergence of HINTS trained with error-based and residual-based losses across different norms for the diffusion and Helmholtz equations.

In both diffusion and Helmholtz equations, residual-based losses outperform error-based losses for HINTS, and the residual-based $\ell_2$ loss achieves the smallest residual level.

For Poisson-like elliptic operators, the residual $r=Ae$ applies a differential operator to the error, which tends to amplify higher-frequency components. As a result, residual-based objectives are often more sensitive to oscillatory errors than error-based objectives. For the indefinite Helmholtz operator, residual-based objectives often amplify error modes away from the resonant band $|\kappa|\approx k$.

Note that for indefinite Helmholtz problems, residual norms may underweight near-resonant components; the observed improved performance is empirical and problem-dependent. Moreover, the preferred training objective is highly architecture-dependent: FNS exhibits a distinct preference for error-based loss functions when trained in a static framework.

\begin{figure}[htbp]
    \centering
    \includegraphics[width=0.60\linewidth]{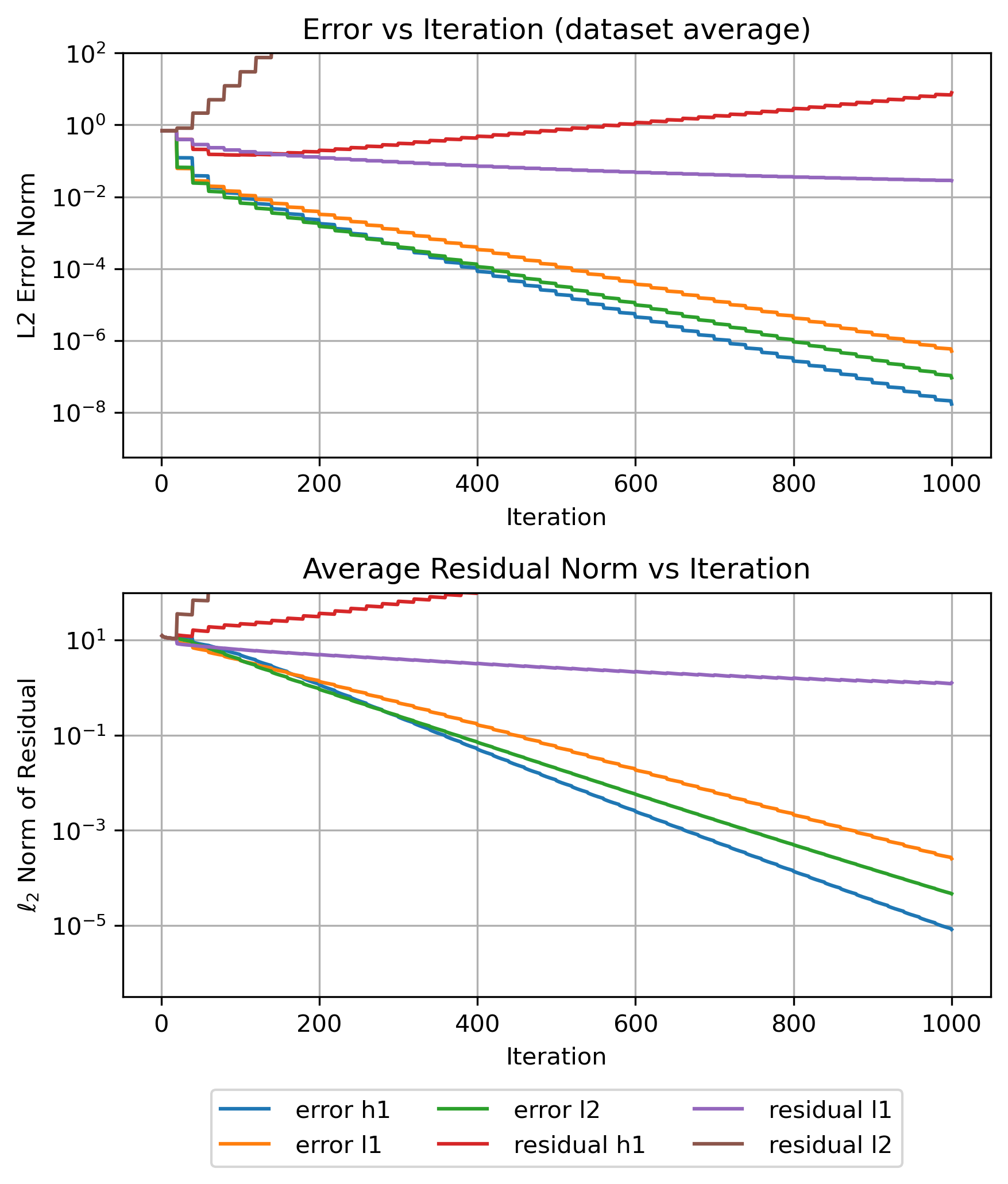}
    \caption{Convergence Behavior of FNS Statically Trained by Residual-Based and Error-Based Losses on Diffusion Equations ($N=201$).}
    \label{fig:loss_fns_diffusion}
\end{figure}

Figure~\ref{fig:loss_fns_diffusion} reports the convergence behavior of FNS under a static training framework with error-based and residual-based objectives.

In contrast to HINTS, for FNS under the static training framework, we observe that residual-based objectives perform substantially worse and even diverge. A plausible explanation is that FNS was originally designed to be trained under a dynamic training framework, whereas the single-step static training objective forces the model to minimize the residual in a one-shot manner, inducing short-sighted and greedy optimization behavior.

Since FNS operates in the spectral domain to mitigate spectral bias~\cite{cui2025hybrid}, residual-based losses may overweight certain spectral components and encourage aggressive corrections, which may lead to instability when repeatedly applied.

\subsection{Cost-Efficiency Analysis of Training Frameworks}
The dynamic training framework, by unrolling the iterative process, can theoretically mitigate the distribution mismatch between the training and inference stages. Meanwhile, since the inference executes the same iterative process, the training frameworks do not incur additional computational costs during the inference stage. The main trade-off lies in the offline training cost in both computation and memory: when unrolling DL-HIM for $K$ steps, backpropagation needs to store intermediate activations and computational graphs for all $K$ steps. Consequently, the training time per epoch and graphics processing unit (GPU) memory usage typically increase linearly with $K$, i.e., $\mathcal{O}(K)$.

Table~\ref{tab:cost_efficiency} compares static and dynamic training (with $K=10$ and curriculum learning) for both HINTS and FNS trained on 1D diffusion equations ($N_\text{train}=31$). In our experiments on 1D ($N_\text{test} = 601$), dynamic training yields only marginal convergence improvements ($1.16\times$--$1.43\times$ at the $10^{-6}$ threshold), while incurring $10.90\times$ training time and $1.61\times$ memory overhead for HINTS, and $5.96\times$ training time and $8.53\times$ memory overhead for FNS. 

Previous and further experiments (see Section~\nameref{sec:experiment_loss_functions} and Section~\nameref{sec:experiment_update}) indicate that, compared with the costly dynamic unrolling strategy, improving the update scheme and choosing suitable objectives are more likely to alleviate convergence stagnation at a much lower cost.

\begin{table}[htbp]
    \centering
    \caption{Cost-efficiency comparison of dynamic versus static training on 1D diffusion equations ($N=601$). All values are normalized to each model's respective static baseline ($1.00\times$).}
    \label{tab:cost_efficiency}
    \begin{tabular}{@{} l cc c cc @{}}
        \toprule
        & \multicolumn{2}{c}{\textbf{HINTS}} & & \multicolumn{2}{c}{\textbf{FNS}} \\
        \cmidrule{2-3} \cmidrule{5-6}
        \textbf{Metric} & \textbf{Static} & \textbf{Dynamic}$^*$ & & \textbf{Static} & \textbf{Dynamic}$^*$ \\
        \midrule
        \multicolumn{6}{@{}l}{\textit{Computational Overheads}} \\
        \ Training time   & 1.00$\times$ & 10.90$\times$ & & 1.00$\times$ & 5.96$\times$ \\
        \ Peak memory     & 1.00$\times$ & 1.61$\times$ & & 1.00$\times$ & 8.53$\times$ \\
        \midrule
        \multicolumn{6}{@{}l}{\textit{Convergence Speedup$^\dagger$}} \\
        \ Error threshold & 1.00$\times$ & 1.16$\times$ & & 1.00$\times$ & 1.21$\times$ \\
        \ Residual threshold& 1.00$\times$ & 1.17$\times$ & & 1.00$\times$ & 1.43$\times$ \\
        \bottomrule
    \end{tabular}
    
    \vspace{4pt}
    \parbox{\columnwidth}{\footnotesize 
        \textit{$^*$Training online using curriculum learning (horizon $K=10$).}\\
        \textit{$^\dagger$Convergence speedup is measured based on the number of iterations to reach the $10^{-6}$ threshold.}
    }
\end{table}

\subsection{Update Strategy}\label{sec:experiment_update}
We further investigate the impact of update strategies on the convergence of DL-HIM. Using the neural operator model trained dynamically by the error-based $\ell_2$ norm, we compared the following four strategies:
\begin{enumerate}
    \item \textbf{Fixed step size}: DL-HIM with a constant step size $\alpha_k=1$ of the neural operator;

    \item \textbf{Exact line search}: DL-HIM with an adaptive step size $\alpha_k = {p_k^\top r_k}/{p_k^\top Ap_k}$ based on exact line search, following FNS~\cite{cui2025hybrid};

    \item \textbf{Standard AA}: DL-HIM accelerated by standard AA (Algorithm~\ref{alg:aa}) with a memory size $m=10$;

    \item \textbf{Physics-aware AA}: DL-HIM accelerated by Physics-Aware AA (Algorithm~\ref{alg:pa_aa}) with a memory size $m=10$;
\end{enumerate}
In AA and PA-AA, the estimates in each step improve with larger $m$, at the cost of a larger least-squares problem. We chose $m = 10$ as the cut-off point, as it was sufficient for speeding up the convergence of DL-HIMs. However, we leave a comprehensive hyperparameter optimization of $m$ to future work.

\begin{figure*}[htbp]
    \centering
    \begin{subfigure}{0.49\linewidth}
        \centering
        \includegraphics[width=1\linewidth]{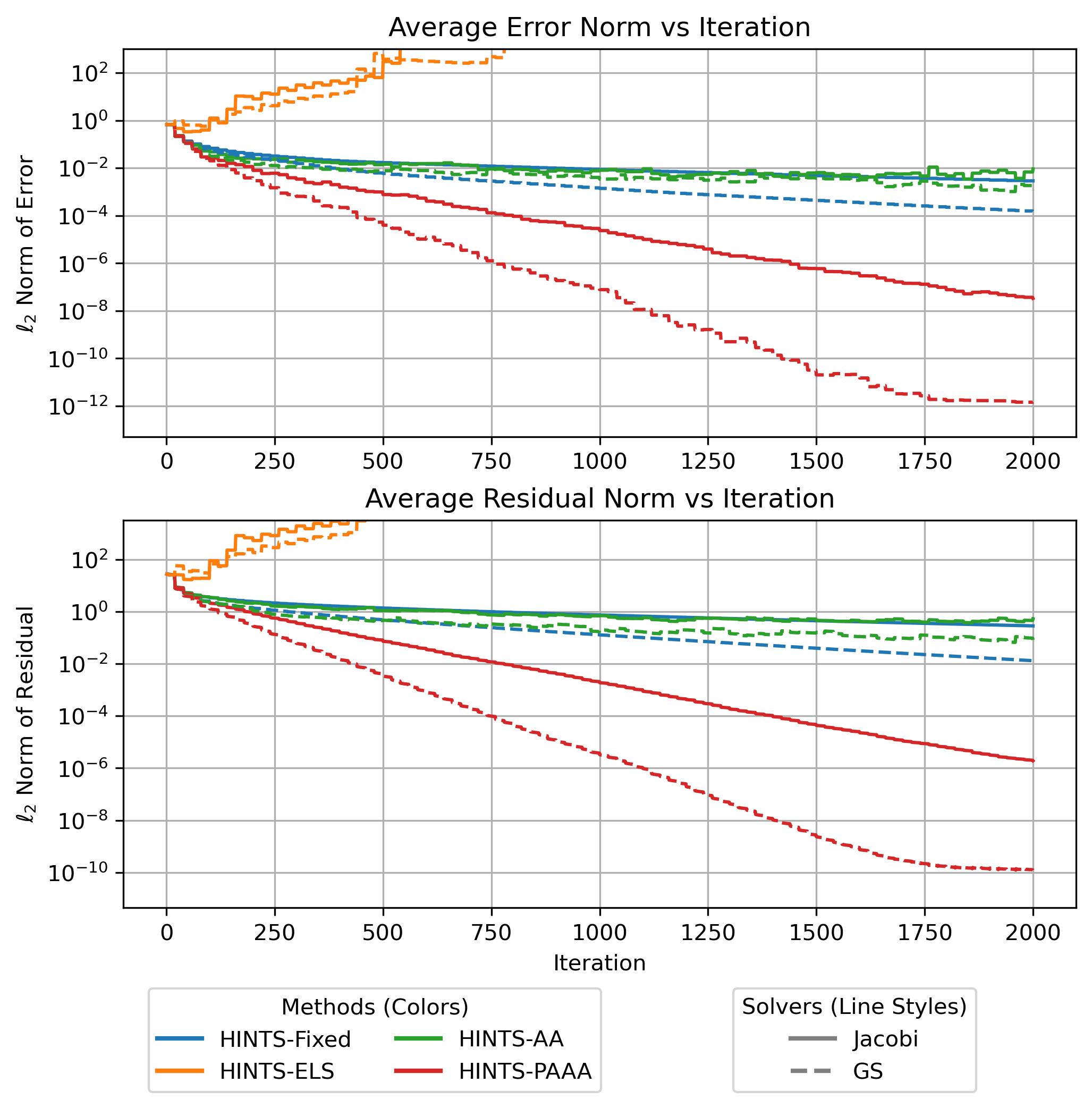}
        \caption{HINTS on 1D Helmholtz ($N=801$)}
    \end{subfigure}
    \hfill
    \begin{subfigure}{0.49\linewidth}
        \centering
        \includegraphics[width=\textwidth]{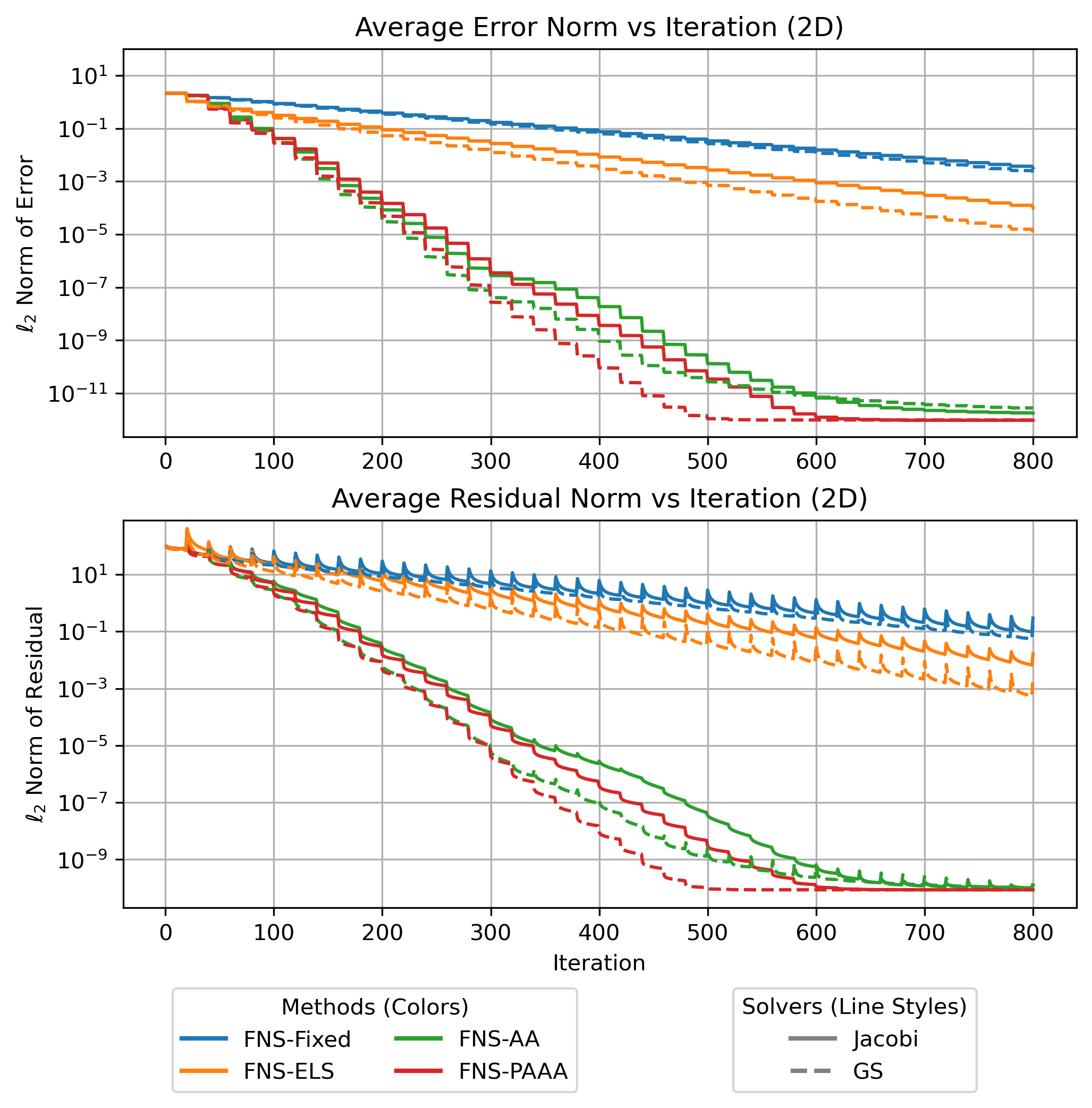}
        \caption{FNS on 2D Diffusion ($N=121{\times}121$)}
    \end{subfigure}
    \caption{Convergence of DL-HIMs with Different Update Strategies: Fixed Step Size (Fixed), Exact Line Search (ELS), Standard AA, and PA-AA for Jacobi and Gauss-Seidel Smoothers}
    \label{fig:aa_performance}
\end{figure*}

The convergence behavior is reported in Figure~\ref{fig:aa_performance} for both Jacobi and Gauss-Seidel smoothers, where models were trained exclusively on independent coarse grids ($N=31$ for 1D and $31\times 31$ for 2D). We observed that the exact line search (ELS) strategy diverges on Helmholtz problems but performs well on Diffusion equations, showing that the adaptive step size is optimal only for SPD quadratic objectives; when applied to non-SPD matrices, this step size is no longer guaranteed to be optimal and may even be counterproductive.

Another critical observation from Figure~\ref{fig:aa_performance} is the performance discrepancy between linear and nonlinear solvers. For HINTS, standard AA does not significantly improve performance. In contrast, FNS preserves linearity with respect to the right-hand side, allowing AA to construct a better linear combination of past updates. As a result, even with a small history size (e.g., $m=10$), AA for FNS still achieves a significant speedup over all other strategies.

By employing PA-AA, which explicitly targets the physical residual, HINTS successfully overcomes this issue. As shown in Figure~\ref{fig:aa_performance}, PA-AA provides a broadly applicable remedy for both nonlinear and linear neural operators. Compared to other update schemes, it significantly accelerates the convergence and enhances the accuracy of DL-HIMs in substantially fewer iterations. The additional performance gain of PA-AA also aligns with the presence of false fixed points in the iterative process of DL-HIMs, which are mitigated when the acceleration is guided by the physical residual rather than the update magnitude.

Finally, we conduct a computational time comparison between the fixed strategy and PA-AA on the 1D Diffusion equations, reported in Table~\ref{tab:speedup}. Employing PA-AA generally speeds up convergence under most error thresholds; this is observed for both the Jacobi and Gauss-Seidel smoothers. Furthermore, although the Gauss-Seidel smoother typically requires fewer iterations to converge than the Jacobi method, it involves a higher computational cost per iteration. Consequently, regarding overall wall-clock computational time, Gauss-Seidel is generally slower than the Jacobi smoother in our experiments.

\begin{table*}
\centering
\caption{Computational time speedup of FNS with Fixed and PA-AA update strategies (Jacobi and Gauss-Seidel) relative to the FNS-Fixed (Jacobi) baseline on 1D diffusion equations($N=1201$), at each representative error and residual thresholds.}
\label{tab:speedup}
\small
\begin{tabular}{clcccc}
\toprule
\textbf{Metric} & \textbf{Threshold} & \thead{FNS-Fixed\\(Jacobi)} & \thead{FNS-PAAA\\(Jacobi)} & \thead{FNS-Fixed\\(GS)} & \thead{FNS-PAAA\\(GS)} \\
\midrule
\multirow{3}{*}{Error}
    & $10^{-2}$ & 1.00$\times$ & 0.61$\times$ & 0.45$\times$ & 0.34$\times$ \\
    & $10^{-4}$ & 1.00$\times$ & 2.43$\times$ & 0.64$\times$ & 1.42$\times$ \\
    & $10^{-6}$ & 1.00$\times$ & 3.25$\times$ & 0.64$\times$ & 1.78$\times$ \\
\midrule
\multirow{3}{*}{Residual}
    & $10^{-2}$ & 1.00$\times$ & 4.14$\times$ & 0.70$\times$ & 2.55$\times$ \\
    & $10^{-4}$ & 1.00$\times$ & 4.41$\times$ & 0.68$\times$ & 2.58$\times$ \\
    & $10^{-6}$ & 1.00$\times$ & 4.50$\times$ & 0.67$\times$ & 2.55$\times$ \\
\bottomrule
\end{tabular}

\vspace{4pt}
\parbox{0.90\textwidth}{\footnotesize\textit{Note: All wall-clock time measurements in this table were conducted using PyTorch 2.6.0 on a single NVIDIA GeForce RTX 4080 GPU.} }

\end{table*}

\section{CONCLUSION}
We have presented a critical evaluation of deep learning-based hybrid iterative methods (DL-HIMs) for solving PDEs, motivated by growing questions about their robustness and reliability in scientific computing. By benchmarking a nonlinear DeepONet-based solver (HINTS) against a linear FFT-based solver (FNS), we showed that convergence failures commonly attributed to neural architectures are also driven by training paradigms and update strategies. Our findings lead to the following practical conclusions.

\paragraph{1. Training Objectives Must Be Matched to Solver Architecture and Physics} The choice of training objective plays a decisive role in convergence behavior. For DeepONet-based solvers, residual-based objectives eliminate the need for expensive reference solutions and consistently lead to more reliable and faster convergence than error-based losses. By contrast, for FFT-based solvers such as FNS, residual-based objectives under static training can be highly sensitive and may even degrade convergence. This highlights that training objectives must be co-designed with the solver architecture and training framework to avoid systematically neglecting intermediate-frequency errors.

\paragraph{2. Dynamic Training Should Be Weighed Against Its Cost}
Although dynamic training frameworks can reduce distribution mismatch by unrolling multiple solver iterations, they introduce substantial computational overhead. In our cost-efficiency study, dynamic training with a horizon $K=10$ incurs $10.90\times$ training and $1.61\times$ memory overhead for HINTS, and $5.96\times$ training time and $8.53\times$ memory overhead for FNS, while yielding only marginal convergence improvements. These results suggest that dynamic training should be adopted with caution. When its benefits are limited, a more effective strategy is to retain static training and focus on lower-cost interventions, such as revisiting the training objective or improving the update mechanism.

\paragraph{3. Robust Update Strategies Are Essential to Avoid Stagnation}
A central failure mode observed across DL-HIMs is convergence to false fixed points, where neural updates vanish while the physical residual remains large. We showed that physics-aware Anderson acceleration (PA-AA), which explicitly minimizes the physical residual rather than the fixed-point update, effectively overcomes this limitation. In our Helmholtz experiments ($N=801$), standard Anderson acceleration (AA) stagnates around $10^{-2}$, whereas PA-AA enables HINTS to converge robustly to $10^{-10}$ in substantially fewer iterations.

\paragraph{Closing Perspective} These results provide a clear answer to current controversies surrounding AI-based PDE solvers. DL-HIMs are not inherently unreliable, but their performance critically depends on physically informed training objectives and update strategies, rather than architectural expressiveness alone. Reliable DL-HIMs therefore require a holistic design in which the network architecture, training paradigm, and iteration strategy are jointly selected based on the PDE, discretization, and computational budget.

\section{ACKNOWLEDGEMENTS}
The authors disclose that an AI-based tool was used for the English improvements in this manuscript. The tool was not used to generate figures, results, or conclusions. All edits were reviewed and approved by the authors, who take full responsibility for the final manuscript.

\def\refname{REFERENCES}
\bibliographystyle{IEEEtran}
\bibliography{ref}

\end{document}